# Recurrence relations for powers of q-Fibonacci polynomials

*Johann Cigler*

Fakultät für Mathematik
Universität Wien
A-1090 Wien, Nordbergstraße 15

johann.cigler@univie.ac.at

**Abstract**
We derive some $q-$analogs of Euler-Cassini-type identities and of recurrence formulas for powers of Fibonacci polynomials.

## 1. Introduction

The Fibonacci numbers $F_n$ are defined by the recurrence relation

$$F_n = F_{n-1} + F_{n-2} \tag{1.1}$$

with initial values $F_0 = 0$ and $F_1 = 1.$

The powers $F_n^k$, $k = 1, 2, 3, \cdots,$ satisfy the recurrence relation

$$\sum_{j=0}^{k+1} (-1)^{\binom{j+1}{2}} \left\langle {k+1 \atop j} \right\rangle F_{n-j}^k = 0, \tag{1.2}$$

where $\left\langle {n \atop k} \right\rangle = \dfrac{\prod_{i=0}^{k-1} F_{n-i}}{\prod_{i=1}^{k} F_i}$ is a so called fibonomial coefficient.

E.g. the squares of the Fibonacci numbers satisfy the recurrence $F_n^2 - 2F_{n-1}^2 - 2F_{n-2}^2 + F_{n-3}^2 = 0.$

The triangle of Fibonomial coefficients ( see A010048 or A055870 in the On-Line Encyclopedia of Integer Sequences [7] ) begins with

```
          1
        1   1
      1   1   1
    1   2   2   1
  1   3   6   3   1
1   5  15  15   5   1
```

The Fibonacci polynomials $F_n(x,s)$ are defined by the recurrence relation

$$F_n(x,s) = xF_{n-1}(x,s) + sF_{n-2}(x,s) \tag{1.3}$$

with initial values $F_0(x,s)=0$ and $F_1(x,s)=1.$
The first terms of this sequence are $0,1,x,x^2+s,x^3+2sx,x^4+3sx^2+s^2,\cdots.$

The powers $F_n^k(x,s)$, $k=1,2,3,\cdots,$ satisfy the recurrence relation

$$\sum_{j=0}^{k+1}(-1)^{\binom{j+1}{2}}s^{\binom{j}{2}}\left\langle {k+1 \atop j} \right\rangle(x,s)F_{n-j}^k(x,s)=0, \qquad (1.4)$$

where the (polynomial-) fibonomial coefficients are defined by

$$\left\langle {n \atop k} \right\rangle(x,s)=\frac{\prod_{i=0}^{k-1}F_{n-i}(x,s)}{\prod_{i=1}^{k}F_i(x,s)}. \qquad (1.5)$$

E.g. for $k=2$ we get the recurrence relation
$F_n(x,s)^2-(x^2+s)F_{n-1}(x,s)^2-s(x^2+s)F_{n-2}(x,s)^2+s^3F_{n-3}(x,s)^2=0.$

The simplest proof of these facts depends on the Binet formula

$$F_n(x,s)=\frac{\alpha^n-\beta^n}{\alpha-\beta}, \qquad (1.6)$$

where

$$\alpha=\frac{x+\sqrt{x^2+4s}}{2},\beta=\frac{x-\sqrt{x^2+4s}}{2}. \qquad (1.7)$$

From (1.6) it is clear that $F_n(x,s)^k$ is a linear combination of $\alpha^{(k-j)n}\beta^{jn}$, $0\le j\le k$. Let $U$ be the shift operator $Uh(n)=h(n-1).$ The sequences $\left(\alpha^{(k-j)n}\beta^{jn}\right)_{n\ge 0}$ satisfy the recurrence relation $\left(1-\alpha^{k-j}\beta^jU\right)\left(\alpha^{(k-j)n}\beta^{jn}\right)=0.$
Since the operators $1-\alpha^{k-j}\beta^jU$ commute we get

$$\left(\prod_{j=0}^{n}(1-\alpha^{k-j}\beta^jU)\right)F_n(x,s)^k=0. \qquad (1.8)$$

As has been observed by L. Carlitz [2] we can now apply the $q-$binomial theorem (cf. e.g. [4])

$$\prod_{j=0}^{n-1}(1-q^jx)=\sum_{k=0}^{n}(-1)^kq^{\binom{k}{2}}\begin{bmatrix} n \\ k \end{bmatrix}x^k. \qquad (1.9)$$

Here $\begin{bmatrix} n \\ k \end{bmatrix}=\begin{bmatrix} n \\ k \end{bmatrix}(q)=\frac{(1-q^n)(1-q^{n-1})\cdots(1-q^{n-k-1})}{(1-q)(1-q^2)\cdots(1-q^k)}$ denotes a $q-$binomial coefficient.

For $q=\dfrac{\beta}{\alpha}$ we get $\begin{bmatrix} n \\ k \end{bmatrix} = \alpha^{k^2-nk}\left\langle \begin{matrix} n \\ k \end{matrix} \right\rangle (x,s).$

This implies

$$\prod_{j=0}^{k}(1-\alpha^{k-j}\beta^{j}U)=\prod_{j=0}^{k}(1-\left(\frac{\beta}{\alpha}\right)^{j}(\alpha^{k}U))=\sum_{j=0}^{k}(-1)^{j}\left(\frac{\beta}{\alpha}\right)^{\binom{j}{2}}\alpha^{j^2-(k+1)j}\left\langle \begin{matrix} k+1 \\ j \end{matrix} \right\rangle (x,s)\alpha^{kj}U^{j}$$

$$=\sum_{j=0}^{k}(-1)^{j}\left(\alpha\beta\right)^{\binom{j}{2}}\left\langle \begin{matrix} k+1 \\ j \end{matrix} \right\rangle (x,s)U^{j}=\sum_{j=0}^{k}(-1)^{\binom{j+1}{2}}s^{\binom{j}{2}}\left\langle \begin{matrix} k+1 \\ j \end{matrix} \right\rangle (x,s)U^{j}.$$

By applying this operator to $F_n(x,s)^k$ we get (1.4).

## 2. Recurrence relations for powers of q-Fibonacci polynomials

The (Carlitz-) $q-$Fibonacci polynomials $f(n,x,s)$ are defined by

$$f(n,x,s)=xf(n-1,x,s)+q^{n-2}sf(n-2,x,s) \tag{2.1}$$

with initial values $f(0,x,s)=0,\ f(1,x,s)=1$ (cf. [3],[5]).

The first values are

$0,1,x,x^2+qs,x^3+(qs+q^2s)x,x^4+(qs+q^2s+q^3s)x^2+q^4s^2,\cdots.$

An explicit expression is

$$f(n,x,s)=\sum_{k\le n-1}\begin{bmatrix} n-1-k \\ k \end{bmatrix}q^{k^2}x^{n-1-2k}s^{k}. \tag{2.2}$$

If we change $q\to\dfrac{1}{q}$ and then $s\to q^{n-1}s$ we get

$f(n-k,x,q^{-k}s)\to f(n-k,x,q^{1-n}s)\to f(n,x,s).$

This implies

### Remark 1

*Each identity*

$$g\left(x,s,q,\ f(n,x,s),\ f(n-1,x,s),\ f(n-2,x,s),\cdots\right)=0 \tag{2.3}$$

*is equivalent with*

$$g\left(x,q^{n-1}s,q^{-1},\ f(n,x,s),\ f(n-1,x,qs),\ f(n-2,x,q^2s),\cdots\right)=0. \tag{2.4}$$

A special case is the well-known fact that (2.1) is equivalent with

$$f(n,x,s) = xf(n-1,x,qs) + qsf(n-2,x,q^2 s). \tag{2.5}$$

The definition of the $q-$Fibonacci polynomials can be extended to all integers such that the recurrence (2.1) remains true. We then get (cf. [5])

$$f(-n,x,s) = (-1)^{n-1} q^{\binom{n+1}{2}} \frac{f(n,x,q^{-n}s)}{s^n}. \tag{2.6}$$

The main aim of this paper is the proof of the following $q-$analog of (1.4) which has been conjectured in [6]:

**Theorem 1**

*Define a $q$-analog of the fibonomial coefficients by*

$$\left\langle {k \atop j} \right\rangle (x,s,q) = \frac{\prod_{i=1}^{k} f(i,x,s)}{\prod_{i=1}^{j} f(i,x,q^{j-i}s) \prod_{i=1}^{k-j} f(i,x,q^j s)}. \tag{2.7}$$

*Then the following recurrence relation holds for all $n \in \mathbb{Z}$ :*

$$\sum_{j=0}^{k+1} (-1)^{\binom{j+1}{2}} s^{\binom{j}{2}} q^{\frac{j(j-1)(2j-1)}{6}} \left\langle {k+1 \atop j} \right\rangle (x,s,q) f(n-j,x,q^j s)^k = 0. \tag{2.8}$$

By Remark 1 this theorem is equivalent to

**Corollary 1**

$$\sum_{j=0}^{k+1} (-1)^{\binom{j+1}{2}} s^{\binom{j}{2}} q^{(n-1)\binom{j}{2} - \frac{j(j-1)(2j-1)}{6}} fibo(k+1,x,s) f(n-j,x,s)^k = 0 \tag{2.9}$$

*with*

$$fibo(k+1,x,s) = \left\langle {k+1 \atop j} \right\rangle (x, q^{n-1}s, q^{-1}) = \frac{\prod_{i=1}^{k} f(i,x,q^{n-i}s)}{\prod_{i=1}^{j} f(i,x,q^{n-j}s) \prod_{i=1}^{k-j} f(i,x,q^{n-i-j}s)}. \tag{2.10}$$

Since there is no $q-$analogue of the Binet formula we use a $q-$analog of an extension of the Cassini - Euler identity for the proof of Theorem 1.

Let

$$fac(k,x,s)=\prod_{i=1}^{k} f(i,x,s) \tag{2.11}$$

and

$$fac(k,x,s,m)=\prod_{i=1}^{k} f(im,x,s). \tag{2.12}$$

Then the following theorem holds:

**Theorem 2**

*For all* $n\in\mathbb{Z}$ *and* $m,\ell\in\mathbb{N}$

$$\det\left(f(n+mi-\ell j,x,q^{\ell j}s)^k\right)_{i,j=0}^{k}=\prod_{j=0}^{k}\binom{k}{j}(-s)^{\binom{k+1}{2}(n-k\ell)+\binom{k+1}{3}(\ell+m)}$$
$$q^{\binom{k+1}{2}\binom{n}{2}+\binom{k+1}{3}mn+\binom{k+1}{3}\frac{m(km-2)}{4}-\binom{k+1}{2}^2\binom{\ell}{2}-\binom{k+1}{3}\frac{\ell(3k+2)}{4}} \tag{2.13}$$
$$\prod_{j=0}^{k-1} fac(k-j,x,q^{mj+n}s,m)\prod_{j=0}^{k-1} fac(k-j,x,q^{\ell j}s,\ell).$$

First we prove the special case $k=1$:

**Lemma 1**

*For all* $n\in\mathbb{Z}$ *and* $m,\ell\in\mathbb{N}$

$$\det\begin{pmatrix} f(n,x,s) & f(n-\ell,x,q^{\ell}s)\\ f(n+m,x,s) & f(n+m-\ell,x,q^{\ell}s)\end{pmatrix}=(-s)^{n-\ell}q^{\binom{n}{2}-\binom{\ell}{2}}f(\ell,x,s)f(m,x,q^n s). \tag{2.14}$$

Various versions of this lemma are well known (cf. [1] and [5]).

Since $f(n+m,x,s)=xf(n+m-1,x,s)+q^{n+m-2}sf(n+m-2,x,s)$ and
$f(n+m-\ell,x,q^{\ell}s)=xf(n+m-1-\ell,x,q^{\ell}s)+q^{n+m-2}sf(n+m-2-\ell,x,q^{\ell}s)$
we see that

$$g(m):=\det\begin{pmatrix} f(n,x,s) & f(n-\ell,x,q^{\ell}s)\\ f(n+m,x,s) & f(n+m-\ell,x,q^{\ell}s)\end{pmatrix}$$

satisfies $g(m)=xg(m-1)+q^{m+n-2}sg(m-2)$ and $g(0)=0.$ Therefore $g(m)=cf(m,x,q^n s)$ for some constant $c.$ To compute $c$ we set $m=-n.$ This gives
$g(-n)=f(n,x,s)f(-\ell,x,q^{\ell}s)=cf(-n,x,q^n s)$ or

$$f(n,x,s)(-1)^{\ell-1}q^{\binom{\ell+1}{2}-\ell^2}s^{-\ell}f(\ell,x,s)=cf(-n,x,q^n s)=c(-1)^{n-1}q^{\binom{n+1}{2}-n^2}s^{-n}f(n,x,s)$$

and therefore

$$g(m)=(-s)^{n-\ell}q^{\binom{n}{2}-\binom{\ell}{2}}f(\ell,x,s)f(m,x,q^n s)=(-s)^{n-\ell}q^{\frac{(n-\ell)(\ell+n-1)}{2}}f(\ell,x,s)f(m,x,q^n s).$$

As a special case we get

**Corollary 2**

*For each* $k \in \mathbb{N}$ *there is a representation of* $f(n-k,x,q^k s)$ *as a linear combination of* $f(n,x,s)$ *and* $f(n-1,x,qs)$:

$$f(n-k,x,q^k s) = \frac{1}{v(k)}\left( f(k-1,x,qs) f(n,x,s) - f(k,x,s) f(n-1,x,qs)\right) \tag{2.15}$$

with

$$v(k) = (-1)^k q^{\binom{k}{2}} s^{k-1}. \tag{2.16}$$

**Proof of Theorem 2**

Using (2.15) we get

$$\det\left( f(n+mi-\ell j,x,q^{\ell j}s)^k \right)_{i,j=0}^{k}$$

$$= \det\left(\left(v(\ell j)^{-1}\left( f(\ell j-1,x,qs) f(n+mi,x,s) - f(\ell j,x,s) f(n+mi-1,x,qs)\right)\right)^k\right)_{i,j=0}^{k}$$

$$= (-1)^{k\ell\binom{k+1}{2}} \frac{1}{s^{(k\ell-2)\binom{k+1}{2}} q^{\left(\binom{\ell}{2}+\binom{2\ell}{2}+\cdots+\binom{k\ell}{2}\right)k}} \det\left(\left(a_j f(n+mi,x,s) + b_j f(n+mi-1,x,qs)\right)^k\right)_{i,j=0}^{k} \tag{2.17}$$

$$= (-1)^{k\ell\binom{k+1}{2}} s^{-(k\ell-2)\binom{k+1}{2}} q^{-\frac{k^2(k+1)\ell(2k\ell+\ell-3)}{12}} \det\left(\left(a_j f(n+mi,x,s) + b_j f(n+mi-1,x,qs)\right)^k\right)_{i,j=0}^{k}$$

with $a_j = f(\ell j-1,x,qs), b_j = -f(\ell j,x,s).$

Since the determinant is multilinear and alternating we get

$$\det\left(\left(a_j f(n+mi,x,s) + b_j f(n+mi-1,x,qs)\right)^k\right)_{i,j=0}^{k}$$

$$= \det\left(\sum_{h=0}^{k}\binom{k}{h}\left(a_j f(n+mi,x,s)\right)^h \left(b_j f(n+mi-1,x,qs)\right)^{k-h}\right)_{i,j=0}^{k}$$

$$= \prod_{j=0}^{k}\binom{k}{j}\sum_{\pi}\det\left(\left(a_j f(n+mi,x,s)\right)^{\pi(j)}\left(b_j f(n+mi-1,x,qs)\right)^{k-\pi(j)}\right)$$

$$= \prod_{j=0}^{k}\binom{k}{j}\sum_{\pi}\det\left(\left(a_j f(n+mi,x,s)\right)^{\pi(j)}\left(b_j f(n+mi-1,x,qs)\right)^{k-\pi(j)}\right)$$

$$= \prod_{j=0}^{k}\binom{k}{j}\sum_{\pi}\prod_{j=0}^{k} a_j^{\pi(j)} b_j^{k-\pi(j)} \det\left(\left(f(n+mi,x,s)\right)^{\pi(j)}\left(f(n+mi-1,x,qs)\right)^{k-\pi(j)}\right)$$

$$= \prod_{j=0}^{k}\binom{k}{j}\sum_{\pi}\operatorname{sgn}(\pi)\prod_{j=0}^{k} a_j^{\pi(j)} b_j^{k-\pi(j)} \det\left(\left(f(n+mi,x,s)\right)^{j}\left(f(n+mi-1,x,qs)\right)^{k-j}\right)$$

$$= \prod_{j=0}^{k}\binom{k}{j}\det\left(a_i^j b_i^{k-j}\right)\det\left(\left(f(n+mi,x,s)\right)^{j}\left(f(n+mi-1,x,qs)\right)^{k-j}\right).$$

Now we need

## Lemma 2

*For* $m \in \mathbb{N}$

$$D(n,m,s,k) = \det\left( f(n+mi,x,s)^j f(n+mi-1,x,qs)^{k-j} \right)$$
$$= (-1)^{\binom{k+1}{2}n+\binom{k+1}{3}m} s^{\binom{k+1}{2}(n-1)+\binom{k+1}{3}m} q^{\binom{k+1}{2}\binom{n}{2}+\binom{k+1}{3}\left(\binom{m}{2}+mn\right)+\binom{k+1}{4}m^2} \tag{2.18}$$
$$\prod_{j=0}^{k-1} fac(k-j,x,q^{mj+n}s,m).$$

## Proof

Using formula $f(n+mi,x,s) = xf(n+mi-1,x,qs) + qsf(n+mi-2,x,q^2s)$
we get as above

$$D(n,m,s,k) = (-s)^{n\binom{k+1}{2}} q^{\binom{k+1}{2}\binom{n+1}{2}} D(0,m,q^n s,k). \tag{2.19}$$

For

$$D(n,m,s,k) =$$
$$\det\left( f(n+mi-1,x,qs)^{k-j} \left( xf(n+mi-1,x,qs) + qsf(n+mi-2,x,q^2s) \right)^j \right)$$
$$= \det\left( f(n+mi-1,x,qs)^{k-j} \left( qsf(n+mi-2,x,q^2s) \right)^j \right)$$
$$= (qs)^{\binom{k+1}{2}} \det\left( f(n+mi-2,x,q^2s)^j f(n+mi-1,x,qs)^{k-j} \right) = (-qs)^{\binom{k+1}{2}} D(n-1,m,qs,k)$$
$$= (-s)^{n\binom{k+1}{2}} q^{\binom{k+1}{2}\binom{n+1}{2}} D(0,m,q^n s,k).$$

Finally we expand $D(0,m,s,k)$ with respect to the first column and get

$$D(0,m,s,k) = \det\left( f(mi,x,s)^j f(mi-1,x,qs)^{k-j} \right)$$
$$= f(-1,x,qs)^k \det \begin{pmatrix} f(m,x,s) f(m-1,x,qs)^{k-1} & f(m,x,s)^2 f(m-1,x,qs)^{k-2} & \cdots & f(m,x,s)^k \\ f(2m,x,s) f(2m-1,x,qs)^{k-1} & f(2m,x,s)^2 f(2m-1,x,qs)^{k-2} & \cdots & f(2m,x,s)^k \\ \cdots & \cdots & \cdots & \cdots \\ f(km,x,s) f(km-1,x,qs)^{k-1} & f(km,x,s)^2 f(km-1,x,qs)^{k-2} & & f(km,x,s)^k \end{pmatrix}$$
$$= f(-1,x,qs)^k f(m,x,s) f(2m,x,s) \cdots f(km,x,s) D(m,m,s,k-1).$$

Thus we have

$$D(0,m,s,k) = f(-1,x,qs)^k f(m,x,s) f(2m,x,s) \cdots f(km,x,s) D(m,m,s,k-1). \tag{2.20}$$

For $k=1$ we get from (2.14) that

$$D(n,m,s,1)=(-1)^n s^{n-1} q^{\binom{n}{2}} f(m,x,q^n s).$$

Therefore Lemma 2 is true for $k=1$.

The general case follows by using (2.19) and (2.20)

$$D(n,m,s,k)=(-s)^{n\binom{k+1}{2}} q^{\binom{k+1}{2}\binom{n+1}{2}} (q^n s)^{-k} f(m,x,q^n s) f(2m,x,q^n s)\cdots f(km,x,q^n s) D(m,m,q^n s,k-1)$$

$$=(-s)^{n\binom{k+1}{2}} q^{\binom{k+1}{2}\binom{n+1}{2}} (q^n s)^{-k} f(m,x,q^n s) f(2m,x,q^n s)\cdots f(km,x,q^n s)$$

$$(-1)^{\binom{k}{2}m+\binom{k}{3}m} (q^n s)^{\binom{k}{2}(m-1)+\binom{k}{3}m} q^{\binom{k}{2}\binom{m}{2}+\binom{k}{3}\left(\binom{m}{2}+m^2\right)+\binom{k}{4}m^2} \prod_{j=0}^{k-2} fac(k-j-1,x,q^{mj+m+n}s,m)$$

$$=(-1)^{n\binom{k+1}{2}+\binom{k+1}{3}m} s^{(n-1)\binom{k+1}{2}+\binom{k+1}{3}m} q^{\binom{k+1}{2}\binom{n}{2}+nm\binom{k+1}{3}+\binom{k+1}{3}\left(\binom{m}{2}\right)+\binom{k+1}{4}m^2} \prod_{j=0}^{k-1} fac(k-j,x,q^{mj+n}s,m).$$

A special case is

**Lemma 3**

$$D(0,\ell,s,k)=\det\left(a_i^j b_i^{k-j}\right)=\det\left(f(\ell i-1,x,qs)^j\left(-f(\ell i,x,s)\right)^{k-j}\right)_{i,j=0}^{k}$$

$$=(-1)^{\binom{k+1}{3}\ell} s^{-\binom{k+1}{2}+\binom{k+1}{3}\ell} q^{\binom{k+1}{3}\binom{\ell}{2}+\binom{k+1}{4}\ell^2} \tag{2.21}$$

$$\prod_{j=0}^{k-1} fac(k-j,x,q^{\ell j}s,\ell).$$

With the use of these lemmas we get

$$\det\left(f(n+mi-\ell j,x,q^{\ell j}s)^k\right)_{i,j=0}^{k}$$

$$=(-1)^{k\ell\binom{k+1}{2}} s^{-(k\ell-2)\binom{k+1}{2}} q^{-\frac{k^2(k+1)\ell(2k\ell+\ell-3)}{12}} \det\left(\left(a_j f(n+mi,x,s)+b_j f(n+mi-1,x,qs)\right)^k\right)_{i,j=0}^{k}$$

$$=(-1)^{k\ell\binom{k+1}{2}} s^{-(k\ell-2)\binom{k+1}{2}} q^{-\frac{k^2(k+1)\ell(2k\ell+\ell-3)}{12}} \prod_{j=0}^{k}\binom{k}{j} D(n,m,s,k)D(0,\ell,s,k)$$

$$=(-1)^{k\ell\binom{k+1}{2}} s^{-(k\ell-2)\binom{k+1}{2}} q^{-\frac{k^2(k+1)\ell(2k\ell+\ell-3)}{12}} \prod_{j=0}^{k}\binom{k}{j} (-1)^{\binom{k+1}{2}n+\binom{k+1}{3}m} s^{\binom{k+1}{2}(n-1)+\binom{k+1}{3}m} q^{\binom{k+1}{2}\binom{n}{2}+\binom{k+1}{3}\left(\binom{m}{2}+mn\right)+\binom{k+1}{4}m^2}$$

$$\prod_{j=0}^{k-1} fac(k-j,x,q^{mj+n}s,m)(-1)^{\binom{k+1}{3}\ell} s^{-\binom{k+1}{2}+\binom{k+1}{3}\ell} q^{\binom{k+1}{3}\binom{\ell}{2}+\binom{k+1}{4}\ell^2}$$

$$\prod_{j=0}^{k-1} fac(k-j,x,q^{\ell j}s,\ell)$$

$$= \prod_{j=0}^{k} \binom{k}{j} (-s)^{\binom{k+1}{2}(n-k\ell)+\binom{k+1}{3}(\ell+m)} q^{\binom{k+1}{2}\binom{n}{2}+\binom{k+1}{3}mn+\binom{k+1}{3}\frac{m(km-2)}{4}-\binom{k+1}{2}^2\binom{\ell}{2}-\binom{k+1}{3}\frac{\ell(3k+2)}{4}}$$

$$\prod_{j=0}^{k-1} fac(k-j,x,q^{mj+n}s,m) \prod_{j=0}^{k-1} fac(k-j,x,q^{\ell j}s,\ell).$$

Thus Theorem 2 is proved.

We will also need some modifications of these results.
Let

$$d(n,m,s,k,j) = \det\Big( f(n+m(i+[i \geq j]),x,s)^j \, f(n+m(i+[i \geq j])-1,x,qs)^{k-j} \Big), \tag{2.22}$$

where $[P]$ denotes the Iverson symbol, i.e. $[P]=1$ if property $P$ is true and $[P]=0$ else.
Then

$$d(n,m,s,k,0) = D(n+m,m,s,k). \tag{2.23}$$

From (2.18) we get

$$d(0,m,s,k,0) = (-1)^{\binom{k+1}{2}m} s^{\binom{k+1}{2}m-k} q^{\frac{km((km+m+k-3)}{4}} fac(k,x,q^m s,m) d(0,m,q^m s,k-1,0). \tag{2.24}$$

Furthermore we get in the same way as above for $j>0$

$$d(0,m,s,k,j) = s^{-k} \frac{fac(k+1,x,s,m)}{f(jm,x,s)} (-s)^{m\binom{k}{2}} q^{\binom{k}{2}\binom{m+1}{2}} d(0,m,q^m s,k-1,j-1). \tag{2.25}$$

## Proof of Theorem 1

The above argument implies that

$$\det\Big( f(n+i-j,x,q^j s)^k \Big)_{i,j=0}^{k+1} = 0. \tag{2.26}$$

If we denote by $A_j$ the matrix obtained by crossing out the first row and the $j-$th column of $\Big( f(n+i-j,x,q^j s)^k \Big)_{i,j=0}^{k+1}$, we get

$$\sum_{j=0}^{k+1} f(n-j,x,q^j s)^k (-1)^j \det\big(A_j\big) = 0$$

or

$$\sum_{j=0}^{k+1} f(n-j,x,q^j s)^k (-1)^j \frac{\det\big(A_j\big)}{\det\big(A_0\big)} = 0. \tag{2.27}$$

To compute $\det\big(A_j\big)$ we use the same method as in Theorem 1.
We get

$$\det\big(A_j\big) = \left( \frac{v(j)}{v(0)v(1)\cdots v(k+1)} \right)^k \det\Big( \big(a_h(j,s) f(n+i,x,s) + b_h(j,s) f(n+i-1,x,qs)\big)^k \Big)_{i,h=0}^{k},$$

where

$a_h(j,s) = f(h-1,x,qs), b_h(j,s) = -f(h,x,s)$ for $h < j$ and
$a_h(j,s) = f(h,x,qs), b_h(j,s) = -f(h+1,x,s)$ for $h \geq j.$

Therefore

$$\det\left(a_h(j,s)^i b_h(j,s)^{k-i}\right) = d(0,1,s,k,j). \tag{2.28}$$

By (2.24) we have
$d(0,1,s,k,0) = (-1)^{\binom{k+1}{2}} s^{\binom{k}{2}} q^{\binom{k}{2}} fac(k,x,qs) d(0,1,qs,k-1,0).$

For $j > 0$ we get from (2.25)

$$d(0,1,s,k,j) = s^{-k} \frac{fac(k+1,x,s)}{f(j,x,s)} (-s)^{\binom{k}{2}} q^{\binom{k}{2}} d(0,1,qs,k-1,j-1). \tag{2.29}$$

Therefore

$$\frac{d(0,1,s,k,j)}{d(0,1,s,k,0)} = (-s)^{-k} \frac{fac(k+1,x,s)}{f(j,x,s)\, fac(k,x,qs)} \frac{d(0,1,qs,k-1,j-1)}{d(0,1,qs,k-1,0)}. \tag{2.30}$$

This implies
$$\frac{d(0,1,s,k,j)}{d(0,1,s,k,0)} = (-s)^{-\sum_{i=0}^{j-1}(k-i)} q^{-\sum_{i=0}^{j-1} i(k-i)} \left\langle {k+1 \atop j} \right\rangle (x,s,q) \frac{d(0,1,q^j s,k-j,0)}{d(0,1,q^j s,k-j,0)}$$
$$= (-s)^{-kj+\binom{j}{2}} q^{-k\binom{j}{2}+\sum_{i=0}^{j-1} i^2} \left\langle {k+1 \atop j} \right\rangle (x,s,q).$$

Therefore we get
$$(-1)^j \frac{\det\left(A_j\right)}{\det\left(A_0\right)} = (-1)^j \left(\frac{v(j)}{v(0)}\right)^k \frac{d(0,1,s,k,j)}{d(0,1,s,k,0)} = (-1)^{j+kj} s^{kj} q^{k\binom{j}{2}} (-s)^{-kj+\binom{j}{2}} q^{-k\binom{j}{2}+\sum_{i=0}^{j-1} i^2} \left\langle {k+1 \atop j} \right\rangle (x,s,q)$$
$$= (-1)^{\binom{j+1}{2}} s^{\binom{j}{2}} q^{\sum_{i=0}^{j-1} i^2} \left\langle {k+1 \atop j} \right\rangle (x,s,q).$$

Thus Theorem 1 is proved.

If we use the fact that $\det\left( f(n+\ell(i-j),x,q^{j\ell}s)^k \right)_{i,j=0}^{k+1} = 0,$ we get in the same way that

$$\sum_{j=0}^{k+1} f(n-j\ell,x,q^{j\ell}s)^k (-1)^j \frac{\det\left(B_j\right)}{\det\left(B_0\right)} = 0,$$

where we denote by $B_j$ the matrix obtained by crossing out the first row and the $j-$th column of $\left( f(n+\ell(i-j),x,q^{j\ell}s)^k \right)_{i,j=0}^{k+1}.$

Here we have

$$(-1)^j \frac{\det\left(B_j\right)}{\det\left(B_0\right)} = (-1)^j \left( \frac{v(j\ell)}{v(0)} \right)^k \frac{d(0,\ell,s,k,j)}{d(0,\ell,s,k,0)} = (-1)^{j+kj\ell} s^{kj\ell} q^{\binom{j\ell}{2}k} \frac{d(0,\ell,s,k,j)}{d(0,\ell,s,k,0)}.$$

If we define

$$\left\langle {k \atop j} \right\rangle (x,s,q,\ell) = \frac{\prod_{i=1}^{k} f(i\ell,x,s)}{\prod_{i=1}^{j} f(i\ell,x,q^{(j-i)\ell}s) \prod_{i=1}^{k-j} f(i\ell,x,q^{j\ell}s)} \tag{2.31}$$

we get from (2.24) and (2.25)

$$\frac{d(0,\ell,s,k,j)}{d(0,\ell,s,k,0)} = \frac{s^{-k} \dfrac{fac(k+1,x,s,\ell)}{f(j\ell,x,s)} (-s)^{\ell\binom{k}{2}} q^{\binom{k}{2}\binom{\ell+1}{2}} d(0,\ell,q^{\ell}s,k-1,j-1)}{(-1)^{\binom{k+1}{2}\ell} s^{\binom{k+1}{2}\ell-k} q^{\frac{k\ell((k\ell+\ell+k-3)}{4}} fac(k,x,q^{\ell}s,\ell) d(0,\ell,q^{\ell}s,k-1,0)}$$

$$= \frac{(-1)^{k\ell}}{s^{k\ell}} q^{-k\binom{\ell}{2}} \frac{fac(k+1,x,s,\ell)}{f(j\ell,x,s)\, fac(k,x,q^{\ell}s,\ell)} \frac{d(0,\ell,q^{\ell}s,k-1,j-1)}{d(0,\ell,q^{\ell}s,k-1,0)}$$

$$= (-1)^{kj\ell-\ell\binom{j}{2}} s^{\ell\binom{j}{2}-kj\ell} q^{-\binom{\ell}{2}\left(kj-\binom{j}{2}\right)-\ell^2\sum_{i=0}^{j-1} i(k-i)} \left\langle {k+1 \atop j} \right\rangle (x,s,q,\ell).$$

Therefore

$$(-1)^j \frac{\det\left(B_j\right)}{\det\left(B_0\right)} = (-1)^{j+kj\ell} s^{kj\ell} q^{\binom{j\ell}{2}k} (-1)^{kj\ell-\ell\binom{j}{2}} s^{\ell\binom{j}{2}-kj\ell} q^{-\binom{\ell}{2}\left(kj-\binom{j}{2}\right)-\ell^2\sum_{i=0}^{j-1} i(k-i)} \left\langle {k+1 \atop j} \right\rangle (x,s,q,\ell)$$

$$= (-1)^{j+\ell\binom{j}{2}} (q^{\frac{(4j+1)\ell-3}{6}} s)^{\ell\binom{j}{2}} \left\langle {k+1 \atop j} \right\rangle (x,s,q,\ell).$$

Thus we get

**Theorem 3**

*For* $k,\ell\in\mathbb{N}$ *the following recurrence relation holds:*

$$\sum_{j=0}^{k+1}(-1)^{j+\ell\binom{j}{2}}\left(q^{\frac{(4j+1)\ell-3}{6}}s\right)^{\ell\binom{j}{2}}\left\langle{k+1 \atop j}\right\rangle(x,s,q,\ell)\, f(n-j\ell,x,q^{j\ell}s)^{k}=0. \qquad (2.32)$$

For the special case $k=1$ this reduces to

$$f(n,x,s)-\frac{f(2\ell,x,s)}{f(\ell,x,q^{\ell}s)}f(n-\ell,x,q^{\ell}s)+(-1)^{\ell}q^{\frac{\ell(3\ell-1)}{2}}s^{\ell}\frac{f(\ell,x,s)}{f(\ell,x,q^{\ell}s)}f(n-2\ell,x,q^{2\ell}s)=0, \qquad (2.33)$$

which has already been proved in [6].